\documentclass[a4paper, 10pt]{article}

   \usepackage[intlimits]{amsmath}
   \usepackage{amssymb}
   \usepackage{bbm}
   \usepackage{amsthm}

\usepackage{amssymb}
\usepackage{amsmath}
\usepackage{graphicx, float, epsfig, url}
\usepackage[english]{babel}
\usepackage[T1]{fontenc}
\usepackage{graphics, enumerate}

\newtheorem{theorem}{Theorem}[section]
\newtheorem{lemma}[theorem]{Lemma}
\newtheorem{proposition}[theorem]{Proposition}
\newtheorem{corollary}[theorem]{Corollary}
\newtheorem{e-definition}[theorem]{Definition\rm}
\newtheorem{remark}{\it Remark\/}
\newtheorem{example}{\it Example\/}

\begin{document}


\author{Andami Ovono Armel}
\title{Asymptotic behaviour for a diffusion equation governed by nonlocal interactions}




\maketitle

\begin{abstract}
In this paper we study the asymptotic behaviour of a nonlocal nonlinear parabolic equation governed by a parameter. After giving the  existence of unique branch of solutions composed by stable solutions in stationary case, we gives for the parabolic problem  $L^\infty $ estimates  of solution based on using the Moser iterations and existence of global attractor. We finish our study by the issue of asymptotic behaviour in some cases when $t\rightarrow \infty$.

\end{abstract}
\section{Introduction}




The non-local issues are important in studying the behavior of certain physical phenomena and population dynamics. A major difficulty in studying these problems often lie in the absence of well-known properties as maximum principle, regularity and properties of Lyapunov (see \cite{ref13}, \cite{ref4}) and also the difficulty to characterize and determine the stationary solutions associated thus making study the asymptotic behavior of these solutions very difficult.$\\$
\noindent In this paper we study the solution $u(t,x)$ to the nonlocal equation
\begin{equation}
\label{eq0}           
\left\{
\begin{aligned}
 u_t -div(a(l_r( u(t)))\nabla u)=f &\quad{dans}\quad \mathbb{R}^+\times\Omega\\
  u(x,t)=0 &\quad{sur}\quad\mathbb{R}^+\times\partial\Omega\\
  u(.,0)=u_0 & \quad{dans}\quad\Omega.
\end{aligned}
\right.
\end{equation}In the above problem $u_0$ and $f$ are such that 
\begin{equation}
u_0\in L^2(\Omega),\quad f\in  L^2(0,T, L^2(\Omega)),
\end{equation}with $T$ a arbitrary positive number, $a$ is a continuous function such that
\begin{equation}
\label{p1}
\exists m,M\quad\text{such that}\quad 0<m\leq a(\epsilon)\leq M\quad \forall\epsilon\in \mathbb{R}.
\end{equation}The nonlocal functional $l_r$ is defined such that
\begin{equation}
l_r(.)(x): L^2(\Omega)\rightarrow\mathbb{R},\quad u\rightarrow l_r(u(t))(x)=\int_{\Omega\cap B(x,r)} g(y)u(t,y) dy.
\end{equation}Here  $ B(x,r)$ is the closed ball of $\mathbb{R}^n$ with radius $r$ and $g\in L^2(\Omega).$.
It is sometimes possible to consider $g$ more generally, especially when one is interested in the study of stationary solutions of (see\,\cite{ref3}).

In physical point of view problem\,(\ref{eq0}) gives many applications especially where $g=1$ in population dynamics. Indeed, in this situation $u$ may represent a population density and $l_r(u)$ the total mass of the subdomain $\Omega\cap B(x,r)$ of  $\Omega$. Hence (\ref{eq0}) can describe the evolution of a population whose diffusion velocity depends on the total mass of a subdomain of $\Omega$. For more details of modelisation we refer the reader to \cite{ref5}. 
This type of equation can be applied more generally to other models including the study of propagation of mutant gene (see \cite{ref8},\cite{ref9},\cite{ref10}). A very recent study of this propagation was made by Bendahmane and Sep{\'u}lveda \cite{ref16} in which they analyze using a finite volume scheme adapted, the transmission of this gene through 3 types of people: susceptible, infected and recovered.

In mathematical point of view, when $r=d$ where $d$ is the diameter of $\Omega$ problem\,(\ref{eq0}) has been studied in various forms (see \cite{ref4},\cite{ref6},\cite{ref7},\cite{ref11}).

However when $0<r<d$, several questions from the theory of bifurcations have arisen concerning the structure of stationary solutions including the existence of a principle of comparison of different solutions depending on the parameter $r$ and the existence of branches (local and global) of solutions. A large majority of these issues has been resolved in \cite{ref3}. It shows that when $a$ is decreasing the existence of a unique global branch of solutions and existence of branch of solutions that are purely local. Some questions may then arise:
\begin{enumerate}[(i)]
\item  The unique branch described in \cite{ref3} it is composed of stable solutions?
\item What about stability properties of the corresponding parabolic problem?
\end{enumerate}The plan for this work is the following. In section\,2 we give some existence and uniqueness results.  Section\,3 is devoted to stationary problem corresponding to (\ref{eq0}). In particular, we study in a radial case, a generalisation of Chipot-Lovat results about determination of the number of solutions. We also establish that the unique global branch of solutions described in \cite{ref3} is composed by stable solutions (theorem\,\ref{bab2}). In section\,4 firstly we address a $L^{\infty}$ estimate taking to account $L^p$ estimate based on Moser iterations. Secondly we prove existence of absorbing set in $H^1_0$, which allows us to prove the existence of a global attractor associated to (\ref{eq0}) (see remark\,\ref{bab}). Finally we obtain a result of stability properties of the corresponding parabolic problem.  
\section{Existence and uniqueness results}
In this section we show a result of existence . We set $ V= H_{0}^{1}(\Omega)$ and  $V'$ its dual, we take the norm in $V$, $\|.\|_V$ such that 
   $$\|u\|^2_V=\int_{\Omega}|\nabla u|^2dx$$
$<.,.>$  means the duality bracket of $V'$ and  $ V.\\$

\noindent Then we have

\begin{theorem}
\label{terib1}  
Let $T>0$,  $f\in  L^2(0,T,V')$ and $u_0\in L^2(\Omega)$, we assume that $a$ is a continuous function and the assumption (\ref{p1}) checked then for every $r$ fixed,  $r\in[0,diam(\Omega)]$, there exists a function $u$ such that
   \begin{equation}
    \label{eq5}
   \begin{cases}
           u\in L^2(0,T,V),\qquad u_t\in L^2(0,T,V')\\
             u(0,.)=u_0\qquad in\qquad\Omega \\
  \frac{d}{dt}(u,\phi)+\int_{\Omega} a(l_{r}( u(t)))\nabla u\nabla\phi dx=<f,\phi>\qquad in\quad D'(0,T)\quad \forall\phi\in H_{0}^{1}(\Omega). 
   \end{cases}
 \end{equation}
Moreover if  $a$ is locally Lipschitz i.e 
 \begin{equation}
 \label{eq6}    
 \forall c \quad\exists \gamma_c\quad\text{such that}\quad|a(\epsilon)-a(\epsilon^\prime)|\leq \gamma_c|\epsilon-\epsilon^\prime|\qquad\forall\epsilon,\,\epsilon^\prime\in [-c,c],
\end{equation} then the solution of (\ref{eq5}) is unique.
\end{theorem}

\begin{remark}
Before to do the proof, it is necessary to see that for  $r=0$ problem (\ref{eq5}) is linear and the proof follows a well-known result (see \cite{ref15}), it is even when  $r=diam(\Omega)$ (see \cite{ref5}). We will focus therefore in the following where  $r\in]0,diam(\Omega)[$.
\end{remark}

\begin{proof}For the existence proof we will use the Schauder  fixed point theorem. Let $w\in L^2(0,T,L^2(\Omega))$ we get  $$t\longrightarrow l_r(w(t)),$$ is measurable as  $a$ is continuous then $$t\longrightarrow a( l_r(w(t))),$$ is too. The problem of finding  $u=u(t,x)$ solution of\begin{equation}
 \label{eq7}
 \begin{cases}
           u\in L^2(0,T,V)\cap C([0,T],L^2(\Omega))\qquad u_t\in L^2(0,T,V')\\
             u(0,.)=u_0 \\
  \frac{d}{dt}(u,\phi)+\int_{\Omega} a(l_{r}( w(t)))\nabla u\nabla\phi dx=<f,\phi>\qquad in\quad D'(0,T)\quad \forall\phi\in H_{0}^{1}(\Omega), 
 \end{cases}
 \end{equation}
is linear, besides (\ref{eq7}) admits a unique solution $u=F_r(w)$  (see\cite{ref15},\cite{ref5}). Thus we show that the application  
   \begin{equation}
   \label{eq8}
   w\longrightarrow F_r(w)=u,
  \end{equation}admits a fixed point. Taking $w=u$ in (\ref{eq7}) we get using  (\ref{p1}) and the Cauchy-Schwarz inequality\begin{equation}
     \label{eq9}
\frac{1}{2}\frac{d}{dt}|u|^2_2+m\| u\|_V^2\leq|f|_{\star}\|u\|_V,
    \end{equation}  $ \|.\|_V$ is the usual norm in $V$ and  $ |f|_{\star}$ is the dual norm of $f.$ We take   $$|u|_{L^2(0,T,V)}=\bigg\{\int_{0}^T\|u\|^2_V dt\bigg\}^\frac{1}{2}.$$Using Young's inequality to the right-hand side of (\ref{eq9}), it follows that\begin{equation}
\label{eq10}
       \frac{1}{2}\frac{d}{dt}|u|^2_2+\frac{m}{2}\| u\|_V^2\leq\frac{1}{2m}|f|_{\star}^2.
\end{equation} 
By integrating  (\ref{eq10}) on $(0,t)$ for $t\leq T$we obtain 
\begin{equation}
\label{eq11} 
    \frac{1}{2}|u(t)|^2_2+\frac{m}{2}\int_{0}^t\| u\|_V^2dt\leq\frac{1}{2}|u_0|^2_2+\frac{1}{2m}\int_{0}^t|f|_{\star}^2.
\end{equation} We deduce that there exists a constant  $ C=C(m,u_0,f)$ such that 
\begin{equation}
\label{eq12}
|u|_{L^2(0,T,V)}\leq C
\end{equation}Moreover 
$$ < u_t,v>+< -div(a(l_{r}( u(t)))\nabla u),v>=<f,v>\quad \forall v\in V,$$ This gives us 
\begin{equation}
\label{eq13}
 |u_t|_{\star}\leq M\|u\|_V+|f|_{\star}.
\end{equation}By raising  (\ref{eq13}) squared and using the Young inequality we have that 
\begin{equation}
\label{eq14}
  |u_t|_{\star}^2\leq 2M^2\|u\|_V^2+2|f|_{\star}^2.
\end{equation}
By integrating (\ref{eq14}) on $ (0,t)$ and assuming (\ref{eq12}) we obtain  \begin{equation}
\label{eq15}
|u_t|_{L^2(0,T,V')}\leq C',
\end{equation}with  $C'=C'(m,M,f,u_0)$ and $C'$ is independent to $w$. It follows from (\ref{eq12}) and (\ref{eq15})  \begin{equation}
\label{eq16} 
|u_t|^2_{L^2(0,T,V')}+|u|^2_{L^2(0,T,V)}\leq R,  
\end{equation}  with $R=C^2+C'^2.$ From  (\ref{eq12}) and the Poincar\'e inequality it follows that 
\begin{equation}
\label{e1q}
|u|_{L^2(0,T,L^2(\Omega))}\leq R',
\end{equation}
By setting\begin{equation}
        \label{eq20}
        R_1=max(R',R),
        \end{equation}and associating  (\ref{e1q}) and (\ref{eq20}), it follows that the application  $F$ maps the ball  $B(0,R_1)$ of  $L^2(0,T,L^2(\Omega))$ into itself. Moreover the balls of  $ H^1(0,T,V,V')$ are relatively compact in $L^2(0,T,L^2(\Omega))$ (see \cite{ref15} for more details), (\ref{eq16}) clearly shows us that  $F(B(0,R_1)$ is relatively compact in  $B(0,R_1)$ with $$B(0,R_1)=\{u\in L^2(0,T,L^2(\Omega));\quad |u|_{L^2(0,T,L^2(\Omega))}\leq R_1 \}.$$ In order to apply the  Schauder fixed point theorem, as announced, we just need to show that $F$ is continuous from $B(0,R_1)$ to itself. This is actually the case and completes the proof of existence.

We will now discuss the uniqueness assuming of course that assumption (\ref{eq6}) be verified.  Consider  $u_1$ and $u_2$ two solutions (\ref{eq5}), by subtracting one obtains in  $\mathit{D}'(0,T)$
\begin{equation}
\label{eq32}
\frac{d}{dt}(u_1-u_2,v)+ \int_{\Omega}( a(l_{r}( u_1(t))\nabla u_1(t)-a(l_{r}( u_2(t)))\nabla u_2(t))\nabla\phi dx=0\qquad\forall\phi\in H_{0}^{1}(\Omega).  
 \end{equation}Since
\begin{equation}
\begin{split}
 a(l_{r}( u_1(t)))\nabla u_1-a(l_{r}( u_2(t)))\nabla u_2(t) &=( a(l_{r}( u_1(t)))-a(l_{r}( u_2(t)))\nabla u_1(t)\\
                                                            &+ a(l_{r}( u_2(t)))\nabla(u_1(t)- u_2(t)),
\end{split}
\end{equation}
we get
\begin{equation}
\label{p2}
\begin{split}
\frac{d}{dt}(u_1-u_2,v)+ &\int_{\Omega}a(l_{r}( u_2(t)))\nabla(u_1(t)- u_2(t))\nabla\phi dx\\
                         &=- \int_{\Omega}( a(l_{r}( u_1(t)))-a(l_{r}( u_2(t)))\nabla u_1\nabla\phi dx\quad\forall\phi\in H_{0}^{1}(\Omega).
\end{split}
\end{equation}
Moreover $u_1,u_2\in C([0,T],L^2(\Omega))$ there exist $z>0$ such that
   \begin{equation}
     \label{eq33}
l_r(u_1(t)),l_r(u_2(t))\in[-z,z].
\end{equation} Taking $v=u_1-u_2$ in (\ref{p2}), it comes easily by Cauchy-Schwartz inequality and (\ref{eq6})  
\begin{equation}
\label{eq35}
\frac{1}{2}\frac{d}{dt}|u_1-u_2|^2_2+m\|u_1-u_2\|^2_V\leq \gamma|l_{r}( u_1(t))-l_{r}( u_2(t))|\| u_1\|_V\|u_1-u_2\|_V .
\end{equation}
We get in \cite{ref3}  
\begin{equation}
\label{eq36}|l_{r}( u(t))\leq C|B(x,r)\cap\Omega|^{\frac{1}{n\vee 3}}| g|_2| u(t)|_2\leq |\Omega|^{\frac{1}{n\vee 3}}| g|_2| u(t)|_2,
\end{equation}
 where   $C$ a constant, $|\Omega|$  represents the measure of $\Omega$ and  $n\vee 3$ the maximum between the dimension $n$ of $\Omega$ and 3. By using (\ref{eq36}), (\ref{eq35}) and the Young inequality  $$ab\leq\frac{1}{2m}b^2+\frac{m}{2}a^2.$$We deduce
\begin{equation}
\label{eq39}
\frac{d}{dt}|u_1-u_2|^2_2+m\|u_1-u_2\|^2_V\leq p(t)| u_1-u_2|_2^2, 
\end{equation}with $$p(t)=\frac{1}{m} (\gamma C\,|\Omega|^{\frac{1}{n\vee 3}}| g|_2\,\| u_1\|_V\,)^2\in L^1(0,T), $$which leads to
\begin{equation}
\label{eq40}
\frac{d}{dt}|u_1-u_2|^2_2\leq\ p(t)| u_1-u_2|_2^2.
\end{equation}
Multiplying (\ref{eq40}) by  $e^{-\int^t_0 p(s)ds}$ it follows that 
\begin{equation}
\label{eq41}
e^{-\int^t_0 p(s)ds}\frac{d}{dt}|u_1-u_2|^2_2- p(t)e^{-\int^t_0 p(s)ds}| u_1-u_2|_2^2\leq 0. 
\end{equation}Hence
\begin{equation}
\label{eq43}
\frac{d}{dt}\{e^{-\int^t_0 p(s)ds}|u_1-u_2|^2_2\}\leq 0. 
\end{equation} 
This shows that  $t\longmapsto e^{-\int\limits^t_0 p(s)ds}|u_1-u_2|^2_2$ is nonincreasing. Since for $t=0$, $$u_1(0,.)=u_2(0,.)=u_0.$$ This function vanishes at $0$  and nonnegative, we conclude that it is identically zero. This concludes the proof.
\end{proof}

\section{Stationary solutions}
Consider the weak formulation to the stationary problem associated to (\ref{eq0}) 
\begin{equation}           
       (P_r) \begin{cases}
            -div(a(l_{r}( u))\nabla u)=f &\mathrm{dans}\quad \Omega\\
               u\in H^1_0(\Omega).
            \end{cases}
            \end{equation}
\subsection{The case $r=d$} By taking $\phi$ the weak solution of the problem 
\begin{equation}
\begin{cases}
 -\Delta\phi=f &\mathrm{dans}\quad \Omega\\
              \phi\in H^1_0(\Omega),
\end{cases}
\end{equation}  we get due to a Chipot-Lovat \cite{ref6}  results that
\begin{theorem}
\label{chi}
Let $a$ be a mapping from $\mathbb{R}$ into $(0, \infty)$. The problem ($P_d$) has many solutions as the problem in $\mathbb{R}$
\begin{equation}
\label{tot}
\mu\,a(\mu)=l_d(\phi),
\end{equation}with $\mu=l_d(u_d)$.
\end{theorem}

\begin{remark}
Theorem\,\ref{chi} allows us to see where $a$ is increasing that the problem $P_d$ admits a unique solution and determine for a given $a$  the exact number of solutions ($P_d$). However it is difficult or impossible to adapt the proof of the theorem\,\ref{chi} in case $0<r<d.$
\end{remark}
\subsection{The case $0<r<d$}
As announced in the introduction we focus our study to the case of radial solutions of ($P_d$). We will assume $\Omega$ is the open ball of $\mathbb{R}^n$ with radius $d/2$ centered at zero. We set
\begin{equation*}
L^2_r(\Omega)=\{u\in L^2(\Omega)\quad\exists\tilde{u}\in L^2(]0,d/2[)\quad \text{such that}\quad u(x)=\tilde{u}(\|x\|) \},
\end{equation*}and we also assume that 
\begin{equation}
\label{cond}
\begin{split}
& f\in L^2_r(\Omega)\\
& g\in L^2_r(\Omega)\\
& a\in W^{1,\infty}(\mathbb{R}),\,\displaystyle{\inf_{\mathbb{R}}a>0}\\
& f\geq 0\quad\text{a.e in}\quad\Omega\\
& g\geq 0\quad\text{a.e in}\quad\Omega.
\end{split}
\end{equation}

We start by giving in some sense in a linear case a result that will be used later to explain the asymptotic behavior. 
\begin{proposition}
\label{marche}
Let $A,B\in C(\overline{\Omega})$ be positive radial functions such that $A\leq B$ in
$ \overline{\Omega}$ and also $f,h\in  L^2(\Omega)$ two positive radial functions. Let $u\in  H^1_0(\Omega)$ the radial solution to
\begin{equation}
\label{comp1}
-\text{div}(A(x)\nabla u)=f\quad\text{in}\quad\Omega,
\end{equation}and
\begin{equation}
\label{comp2}
-\text{div}(B(x)\nabla u)=h\quad\text{in}\quad\Omega.
\end{equation} Then $f\leq h$ a.e in $\Omega.$
\end{proposition}
\begin{proof}
We proved in \cite{ref3} that if $u$ is a the radial solution of (\ref{comp1}) then for a.e $t$ in $[0,d/2]$,
\begin{equation}
\label{comp3}
\tilde{u}^{\prime}(t)=-\frac{1}{\tilde{A}(t)}\int_0^t (\frac{s}{t})^{n-1}\tilde{f}(s)\,ds.
\end{equation}From (\ref{comp1}), (\ref{comp2}) and (\ref{comp3}) we obtain
\begin{equation*}
\frac{\tilde{B}(t)}{\tilde{A}(t)}\int_0^t (\frac{s}{t})^{n-1}\tilde{f}(s)\,ds=\int_0^t (\frac{s}{t})^{n-1}\tilde{h}(s)\,ds.
\end{equation*}Since  $A\leq B$ in $\overline{\Omega}$ and $f,h\geq 0$ with $f\not\equiv 0,\,h\not\equiv 0$ hence $f\leq g$.

\end{proof}

In a nonlocal case, some results of existence of radial solutions and comparison principle between  $u_r$, $u_d$ and $u_0$ has been demonstrated in \cite{ref3}. It is also proved that if we set for all $r\in[0,d]$
\begin{equation}
   I_r:=[\displaystyle{\inf_{\Omega}\, l_r(\phi)},\displaystyle{\sup_\Omega\, l_r(\phi)}].
\end{equation}Here $\phi$ denotes the solution of 
\begin{equation}
\begin{cases}
 -\Delta\phi=f &\mathrm{dans}\quad \Omega\\
              \phi\in H^1_0(\Omega).
\end{cases}
\end{equation}By the inclusion or not of $I_r $ at an interval of $\mathbb{R}$ we somehow generalize the theorem\,\ref{chi}. 
\begin{lemma}
\label{ler1}
Let  $r\in[0,d]$. Assume that (\ref{cond}) holds true and there exist $0\leq m_1\leq m_2$ such that 
\begin{equation}
   a(m_1)=\displaystyle{\max_{[m_1,m_2]} a}\quad a(m_2)=\displaystyle{\min_{[m_1,m_2]} a}
\end{equation}
\begin{equation}
  I_r\subset [m_1a(m_1),m_2a(m_2)]. 
\end{equation}Then ($P_r$) admits a radial solution $u$ and 
\begin{equation}
m_1\leq l_r(u)\leq m_2\quad\text{a.e in}\quad\Omega.
\end{equation}
\end{lemma}For the proof, we refer the reader to \cite{ref3}.

Generalizing this construction type of the diffusion coefficient $a$ we obtain

\begin{proposition}
\label{babe1}
Let  $r\in[0,d]$. Assume that (\ref{cond}) holds true and there exist an odd integer  $n_1$ and $n_1+1$ positive real numbers  $\{m_i\}_{i=0\ldots n_1}$,\,  with $m_0=0$ and for all $i\in\{0,\ldots,n_1-1\}$ we have $m_{i}<m_{i+1}$. Moreover 
\begin{equation}
\label{con1}
\begin{split}
& a(m_i)=\displaystyle{\max_{[m_i,m_{i+1}]} a};\quad a(m_{i+1})=\displaystyle{\min_{[m_i,m_{i+1}]} a}\quad\forall i\in\{0,2,\ldots,n_1-3,n_1-1\}\\& I_r\subset\displaystyle{\bigcap_{i=0,2,\ldots,n_1-3,n_1-1} [m_ia(m_{i}),m_{i+1}a(m_{i+1})]}\\\end{split}
\end{equation}Then ($P_r$) admits at least $\frac{n_1+1}{2}$ radial solutions $\{u_i\}_{i\in\{0,2,\ldots,n_1-1\}}$ such that 
\begin{equation*}
m_i\leq l_r(u_i)\leq m_{i+1}\quad \forall i\in\{0,2,\ldots,n_1-3,n_1-1\}.
\end{equation*}\end{proposition}
\begin{proof}
The proof here is by induction. Indeed we set
\begin{equation*}
\mathcal{P}_{n_1}=\{\text{ If condition} \, (\ref{con1})\text{ is satisfied then} (P_r)\text{ admits at least} \frac{n_1+1}{2} \text{ solutions.}\}
\end{equation*}
By using lemma\,\ref{ler1} with $m_1=0$ and $m_2=m_1$, it is easy to prove for $n_1=1$ that $\mathcal{P}_{n_1}$ is true. For $n_1>1$, This procedure can be repeated to prove that if  $\mathcal{P}_{n_1-2}$  holds true then $\mathcal{P}_{n_1}$ holds too. 
\end{proof}
 \begin{example}
Let us see a function $a$ satisfying proposition\,\ref{babe1}. For this, we consider the case $n_1=3$ and $r\in(0,d]$. Considering (\ref{cond}) and the strong maximum principle we get  $\text{min}\,I_r>0.$ Taking $$m_1:= 2\,\frac{\text{max}\,I_r}{a(0)},\quad a(m_1):=\frac{a(0)}{2}$$  with $a(0)>0$ and also $a$ decreasing on $[0,m_1]$ then we prove lemma\,\ref{ler1} conditions. 

By repeating this process with $m_2>m_1$ and setting $$a(m_2):= \,\frac{\text{min}\,I_r}{m_2},\quad m_3:=2\,\frac{\text{max}\,I_r}{a(m_2)}$$  with $a(m_3):=\frac{a(m_2)}{2}$ and also $a$ is decreasing on $[m_2,m_3]$. This shows the existence of such $a$.

\begin{figure}[h]
\includegraphics[angle=0,width=10cm,height=8cm]{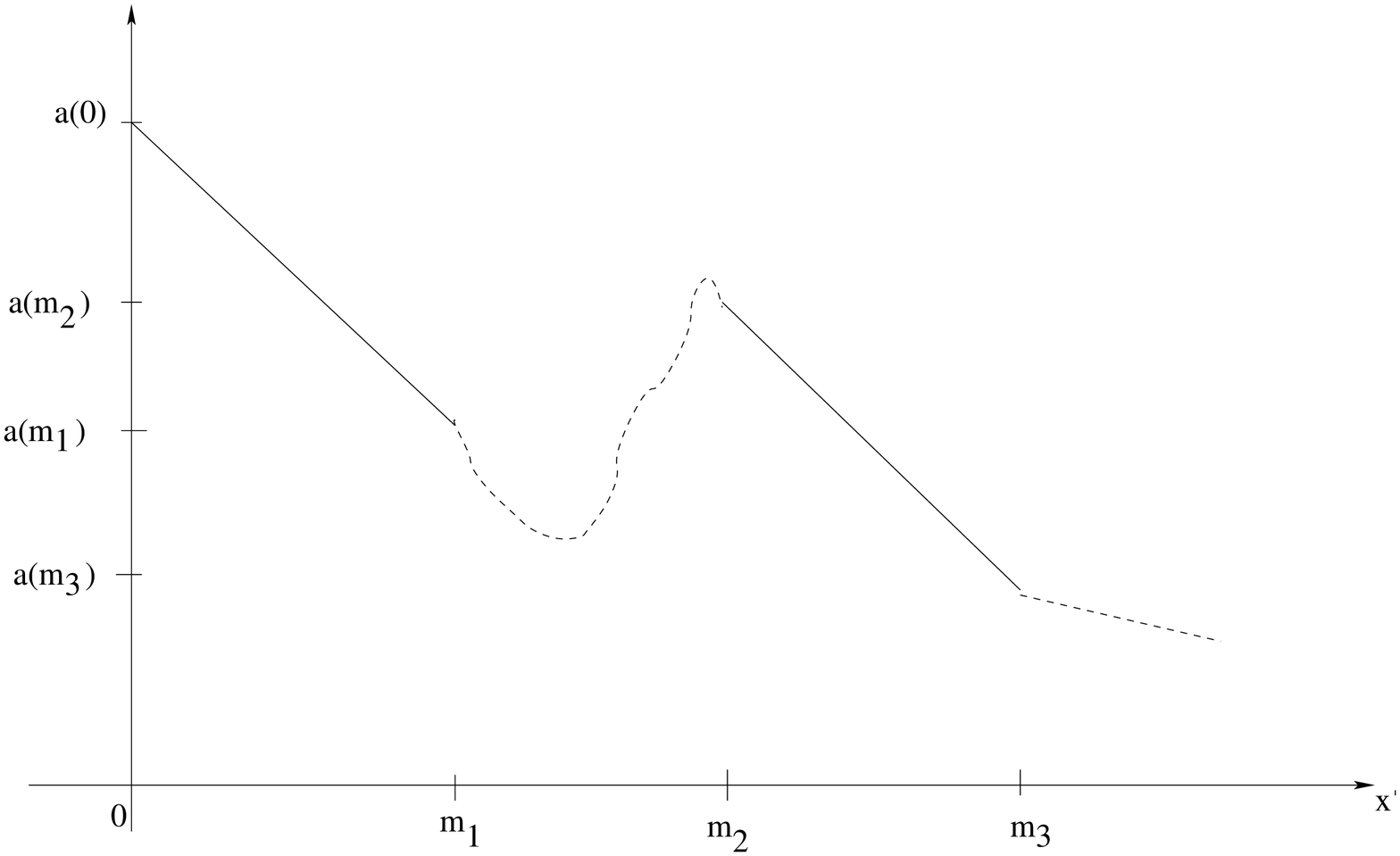}
\caption{The case $n_1=3$}
\label{dor}
\end{figure}

In the representation of $a$ we have deliberately left, on solid line parts of the curve satisfying the conditions of proposition\,\ref{babe1} and dotted line one without constraints. This situation are explain in the figure\,\ref{dor}. 
\end{example}
\begin{remark}
As previously announced, the proposition\,\ref{babe1} generalizes a certain point of view Theorem\,\ref{chi}. However it does not accurately determine the exact number of solutions of ($P_r $) and the bifurcation points of branch of solutions. We have shown in \cite{ref2} way to solve this problem by using the linearized problem, the principle of comparisons obtained in \cite{ref3} and the Krein-Rutman theorem.
\end{remark}
\subsection{Stable solutions of ($P_r$)}
\begin{e-definition}Given  a domain $\Omega\subset\mathbb{R}^n$, a solution $u_r\in H^1_0(\Omega)$   of ($P_r$) is stable if:
\begin{equation}
\label {sta}
\forall\phi\in H^1_0(\Omega)\quad G_{u_r}(\phi):=\int_{\Omega}a(l_r(u_r))|\nabla\phi|^2-\int_{\Omega}a^{\prime}(l_r(u_r))l_r(\phi)\nabla u_r\nabla\phi\geq 0.
\end{equation}
\end{e-definition}
\begin{e-definition}Given $u:\,[0,d]\rightarrow H^1_0(\Omega)$, the graph of $u$ is called a (global) branch of solutions if

\begin{enumerate}[(i)]
\item $u\in C([0,d], H^1_0(\Omega)),$
\item $u(r)$ is solution to ($P_r$) for all $r$ in $[0,d]$.
\end{enumerate}
$u$ is called a local branch if it's defined only on a subinterval of $[0,d]$ with positive measure.
\end{e-definition}
Before concluding this section, we will focus into the case $a$ nonincreasing to prove the stability of the global branch of solutions. 
Assume for all $r\in [0,d]$, $u_r$ is a solution to ($P_r$) and 
\begin{equation}
\label{oli1}
0\leq l_r(u_r)(x)\leq\mu_d\quad\text{for a.e}\quad x\in\Omega.
\end{equation}
Assume that there exists a solution $\mu_d$ to (\ref{tot}) such that
\begin{equation}
\label{tot1}
a(\mu_d)=\displaystyle{\min_{[0,\mu_d]} a}\quad\text{and}\quad a(0)=\displaystyle{\max_{[0,\mu_d]} a}.
\end{equation} We prove in \cite{ref3} 
\begin{theorem}
\label{tata1}
Assume (\ref{cond}), (\ref{oli1}), (\ref{tot1}) and (\ref{tot}) holds. Assume in addition that $a\in W^{1,\infty}(\mathbb{R})$ and for some positive constant $\epsilon$, it holds that
\begin{equation}
\label{const}
C_1|g|_2|f|_2|a^{\prime}|_{\infty,[-\epsilon,\mu_d+\epsilon]}\frac{1}{a(\mu_d)^2}<1,
\end{equation}where $C_1$ is a constant dependent to $\Omega$. Then
\begin{enumerate}[(i)]
\item For all $r$ in $[0,d]$, ($P_r$) possesses a unique radial solution $u_r$ in $[u_0,u_d]$;
\item $\{(r,u_r):r\in[0,d]\}$ is a branch of solutions without bifurcation point;
\item it is only global branch of solutions;
\item if in addition, $a$ is nonincreasing on $[0,\mu_d]$ then $r\mapsto u_r$ is  nondecreasing.

\end{enumerate}
\end{theorem} 
\begin{remark}
It is very difficult to obtain property (iv) for any $a$. However when $a$ is nonincreasing provide us important information for studying the stability of this branch of solutions.
\end{remark}
\begin{corollary}
\label{bab2}Let $u_d^1$ the smallest solution to ($P_d$).
Assume (\ref{cond}) and (\ref{tot}) holds true and there exists a solution $\mu_d$ to (\ref{tot}) satisfied (\ref{tot1}). Assume in addition that $a\in W^{1,\infty}(\mathbb{R})$, $u_d^1$ satisfied (\ref{oli1}) and for some positive constant $\epsilon$, it holds that
\begin{equation}
\label{oublie}
C_1|g|_2|f|_2|a^{\prime}|_{\infty,[-\epsilon,\mu_d+\epsilon]}\frac{1}{a(\mu_d)^2}<1,
\end{equation}where $C_1$ is a constant dependent to $\Omega$. 

Then $\{(r,u_r):r\in[0,d]\}$ is the only global branch of solutions starting to $u_d^1$.

\end{corollary}

\begin{figure}[h]
\includegraphics[angle=0,width=10cm,height=6cm]{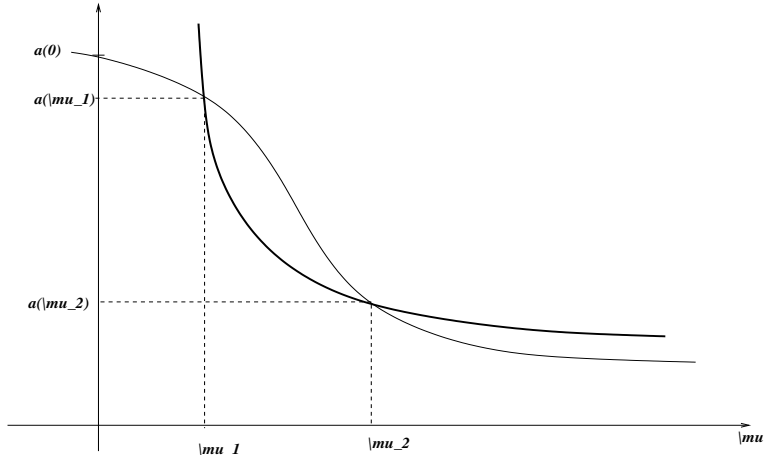}
\caption{case of 2 solutions}
\label{dodo1}
\end{figure}

\begin{proof}The fact that $\{(r,u_r):r\in[0,d]\}$ is the only global branch of solutions results from theorem\,\ref{tata1}.
We will now show that this unique branch of solutions is stable and start at $r=d$ by $u_d^1$. For this we consider without loss of generality ($P_d$) admits two solutions $u_d^1$ and $u_d^2$ such that $u_d^1\leq u_d^2$.
We denote by $\mu_1$ and $\mu_2$ respectively solutions of (\ref{tot}) corresponding to $u_d^1$ and $u_d^2$ (see figure\,\ref{dodo1}). It is easy to see that $\mu_1$ and $\mu_2$ satisfied (\ref{tot1}). 

Assume $\{(r,u_r):r\in[0,d]\}$ is the only global branch of solutions starting to $u_d^2$. Then we get $\quad C_1|g|_2|f|_2|a^{\prime}|_{\infty,[-\epsilon,\mu_2+\epsilon]}\frac{1}{a(\mu_2)^2}<1.$ In this case, using theorem\,\ref{tata1} we get ($P_r$) possesses a unique radial solution $u_r$ in $[u_0,u_d^2]$ and the mapping $r\mapsto u_r$ is nondecreasing. By continuity of this mapping, we can find a $r_0\in ]0,d[$ such that $u_{r_0}=u_d^1\quad\text{for a.e}\quad x\in\Omega$. This means that $u_d^1$ is a solution of ($P_{r_0}$).  This gives us an absurdity and concludes the proof.
\end{proof}

We are now able to prove:

\begin{proposition}
Under assumptions and notation of corollary\,\ref{bab2}, the global branch of solutions described in theorem\,\ref{tata1} is composed by stable solutions.

\end{proposition}
\begin{proof}
For all $r\in[0,d]$, let $u_r$ be a solution belonging to the global branch of solutions described in theorem\,\ref{tata1}. By using the linearized problem of ($P_r$), we get $\forall\phi\in H^1_0(\Omega)$
\begin{equation}
\label{oli2}
\begin{split}
\int_{\Omega}a(l_r(u_r))|\nabla\phi|^2-&\int_{\Omega}a^{\prime}(l_r(u_r))l_r(\phi)\nabla u_r\nabla\phi\geq\\ &\displaystyle{\inf_{\Omega}a(l_r(u_r))}|\nabla\phi|_2^2-C\,|g|_2|a^{\prime}|_{\infty,[-\epsilon,\mu_1+\epsilon]}|\nabla u_r|_2|\nabla\phi|_2^2.
\end{split}
\end{equation}Taking into account that $|\nabla u_r|_2\leq C(\Omega)\frac{|f|_2}{\displaystyle{\inf_{\Omega}a(l_r(u_r))}}$ where  $C(\Omega)$ designed the Poincar\'e Sobolev constant.  We obtain
\begin{equation}
\label{oli3}
\begin{split}
\int_{\Omega}a(l_r(u_r))|\nabla\phi|^2-&\int_{\Omega}a^{\prime}(l_r(u_r))l_r(\phi)\nabla u_r\nabla\phi\geq\\ &|\nabla\phi|_2^2\Big(\displaystyle{\inf_{\Omega}a(l_r(u_r))}-C_1\,|g|_2|a^{\prime}|_{\infty,[-\epsilon,\mu_1+\epsilon]}\frac{|f|_2}{\displaystyle{\inf_{\Omega}a(l_r(u_r))}}\Big).
\end{split}
\end{equation}Moreover by assumptions (\ref{oli1}) and (\ref{tot1}) we get $a(\mu_1)\leq\displaystyle{\inf_{\Omega}a(l_r(u_r))}$. 

Thus (\ref{oublie}) becomes
\begin{equation}
C_1|g|_2|f|_2|a^{\prime}|_{\infty,[-\epsilon,\mu_d+\epsilon]}\frac{1}{\displaystyle{\inf_{\Omega}a(l_r(u_r))}^2}<1.
\end{equation}
We deduces 
\begin{equation}
\int_{\Omega}a(l_r(u_r))|\nabla\phi|^2-\int_{\Omega}a^{\prime}(l_r(u_r))l_r(\phi)\nabla u_r\nabla\phi\geq 0.
\end{equation}This concluded the proof.
\end{proof}

\begin{figure}[h]
\includegraphics[angle=0,width=10cm,height=6cm]{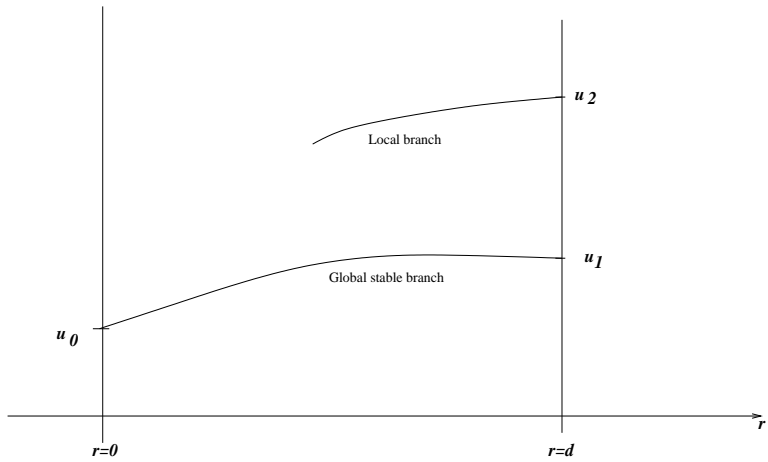}
\caption{Branch of solutions}
\label{dodo2}
\end{figure}

\section{Parabolic problem}
\subsection{$L^{\infty}$  estimate}

In what follows we obtain $L^{\infty}$ estimate of the solution  (\ref{eq0}) from $L^q$ estimate. The method we use is based on iterations Moser type, for more details on the method see \cite{ref12}.

We get
\begin{theorem}
\label{terib3}
Let $n\geq 3$ and $u$ a classical solution of (\ref{eq0}) defined on $[0,T)$. Assume that $p>1$ and $q>1$ such that $\frac{1}{p}+\frac{1}{q}=1$. Suppose further that $U_q=\sup_{t<T}|u(t)|_q<\infty$,  $f\in L^\infty(0,\infty,L^q(\Omega)).$  If $p<\frac{n}{n-2}$ then $U_\infty<\infty$.
\end{theorem}

To prove this theorem we need the following proposition:
\begin{lemma}
\label{moser1}
Consider $u$ a classical solution of (\ref{eq0}) on $[0,T),$\, $r\geq 1$ and $p>1$  such that $\frac{1}{p}+\frac{1}{q}=1$\, with $p<\frac{n}{n-2}$. We take $\tilde U_r=\max\{1,|u_0|_\infty,\, U_r=sup_{t<T}|u(t)|_r\}$ and let $$\sigma(r)=\frac{p(n+2)}{2[r(2p-pn+n)+np]}.$$ Then there exists a constant $C_2=C_2(\Omega,m)$ such that $$\tilde U_{2r}\leq [C_2\,\|f\|_{L^\infty(0,\infty,L^q(\Omega))}]^{\sigma(r)}r^{\sigma(r)}\tilde U_{r}.$$
\end{lemma}
\begin{proof}
Multiplying  (\ref{eq0}) by $u^{2r-1}$ and then using the H\"{o}lder inequality  yields
\begin{equation}
\label{it2}
\frac{1}{2r}\frac{d}{dt}\int_{\Omega}u^{2r}dx+m\frac{2r-1}{r^2}\int_{\Omega}|\nabla (u^r)|^2dx\leq|f|_q|u^{2r-1}|_p.
\end{equation}
As 
\begin{equation}
\label{it3}
|u^{2r-1}|_p=|u^r|^{\frac{2r-1}{r}}_{p\frac{2r-1}{r}},
\end{equation}
by taking $w=u^r$ in (\ref{it2}) and (\ref{it3}), we get easily
\begin{equation}
\label{it4}
\frac{1}{2r}\frac{d}{dt}|w|^2_2+m\frac{2r-1}{r^2}|\nabla w|^2_2\leq|f|_q|w|^{\alpha}_{\alpha p},
\end{equation}with $\alpha=\frac{2r-1}{r}$. Let $\beta$ such that
\begin{equation}
\label{it5}
\frac{1}{\alpha p}=\beta+ \frac{1-\beta}{2^\star},
\end{equation} with $2^\star=\frac{2n}{n-2}.$
We claim that $\beta\in(0,1).\\$
In fact $$\beta=\frac{2nr-(n-2)(2r-1)p}{(n+2)(2r-1)p}.$$ Since $p<\frac{2r}{2r-1}\,\frac{n}{n-2}$ then $\beta>0.$ As well as  $(n+2)(2r-1)p>2nr-(n-2)(2r-1)p$ implies that $\beta<1$ this prove that $\beta\in(0,1).$ 

Using an interpolation inequality (see\,\cite{ref12}) in (\ref{it4}) and (\ref{it5}), we get  
\begin{equation}
\label{it6}
\frac{1}{2r}\frac{d}{dt}|w|^2_2+m\frac{2r-1}{r^2}|\nabla w|^2_2\leq|f|_q\bigg(|w|^{\beta}_1\,|w|^{1-\beta}_{2^\star}\bigg)^{\alpha }.
\end{equation}
Applying Sobolev injections in (\ref{it6}), we have
\begin{equation}
\label{it7}
\frac{1}{2r}\frac{d}{dt}|w|^2_2+m\frac{2r-1}{r^2}|\nabla w|^2_2\leq\bigg[|f|_q\bigg(\frac{2r}{m}\bigg)^{\frac{\alpha(1-\beta)}{2}}|w|^{\beta\alpha}_1\,C^{(1-\beta)\alpha}\bigg]\,\bigg[\bigg(\frac{m}{2r}\bigg)^{\frac{\alpha(1-\beta)}{2}}|\nabla w|^{(1-\beta)\alpha}_2\bigg],\end{equation}
and also
\begin{equation}
\label{it8}
\frac{1}{2r}\frac{d}{dt}|w|^2_2+m\frac{2r-1}{r^2}|\nabla w|^2_2\leq\bigg[|f|_q\bigg(\frac{2r}{m}\bigg)^{\frac{\alpha(1-\beta)}{2}}|w|^{\beta\alpha}_1\,C^{(1-\beta)\alpha}\bigg]\,\bigg[\bigg(\frac{m}{2r}\bigg)|\nabla w|^2_2\bigg]^{\frac{\alpha(1-\beta)}{2}}.
\end{equation}
Since $\beta\in(0,1)$ and $\frac{\alpha}{2}\in (0,1)$ it is clear that $\frac{\alpha(1-\beta)}{2}\in (0,1).$ Applying Young's inequality to (\ref{it8}) with $\frac{\alpha(1-\beta)}{2}+1-\frac{\alpha(1-\beta)}{2}=1$. We obtain
\begin{equation}
\label{it9}
\frac{1}{2r}\frac{d}{dt}|w|^2_2+m\frac{2r-1}{r^2}|\nabla w|^2_2\leq\delta\bigg[|f|_q^{\frac{1}{\delta}}\bigg(\frac{2r}{m}\bigg)^{\frac{\alpha(1-\beta)}{2\delta}}|w|^{\frac{\beta\alpha}{\delta}}_1\,C^{\frac{2}{\delta}}\bigg]+\frac{\alpha(1-\beta)}{2}\bigg[\bigg(\frac{m}{2r}\bigg)|\nabla w|^2_2\bigg],
\end{equation}
with  $\delta=1-\frac{\alpha(1-\beta)}{2}$. 

Joining the fact that $\frac{\alpha(1-\beta)}{2}\in (0,1)$ and $\delta<1$ to  (\ref{it9}), we deduce
\begin{equation}
\label{it10}
\frac{1}{2r}\frac{d}{dt}|w|^2_2+m\frac{3r-2}{2r^2}|\nabla w|^2_2\leq|f|_q^{\frac{1}{\delta}}\bigg(\frac{2r}{m}\bigg)^{\frac{\alpha(1-\beta)}{2\delta}}|w|^{\frac{\beta\alpha}{\delta}}_1\,C^{\frac{2}{\delta}}.
\end{equation}
We set
$$2r\sigma(r)-1=\frac{\alpha(1-\beta)}{2\delta}\quad \text{and}\quad 2\rho(r)=\frac{\beta\alpha}{\delta},$$
 (\ref{it10}) becomes
\begin{equation}
\label{it11}
\frac{1}{2r}\frac{d}{dt}|w|^2_2+m\frac{3r-2}{2r^2}|\nabla w|^2_2\leq|f|_q^{\frac{1}{\delta}}\bigg(\frac{2r}{m}\bigg)^{2r\sigma(r)-1}|w|^{2\rho(r)}_1\,C^{\frac{2}{\delta}}.
\end{equation}
This gives us taking into account that  $\frac{3r-2}{r}>1$

\begin{equation}
\label{it13}
\frac{d}{dt}|w|^2_2+m|\nabla w|^2_2\leq|f|_q^{\frac{1}{\delta}}\bigg(\frac{2r}{m}\bigg)^{2r\sigma(r)}|w|^{2\rho(r)}_1\,m\,C^{\frac{2}{\delta}}.
\end{equation}
By a calculation we can verify that 
\begin{equation*} 
\rho(r)=\frac{2nr-(n-2)(2r-1)p}{2r(p(n+2)+n)-2n(2r-1)p},
\end{equation*}and also that  $\rho(r)\in(0,1).$ 

Using the Poincar\'e Sobolev  inequality and that $\rho(r)<1$ in (\ref{it13}), yields 
\begin{equation}
\label{it14}
\frac{d}{dt}|w|^2_2+\frac{m}{C_1(\Omega)}|w|^2_2\leq|f|_q^{\frac{1}{\delta}}\bigg(\frac{2r}{m}\bigg)^{2r\sigma(r)}|w|^2_1\,m\,C^{\frac{2}{\delta}},
\end{equation}where  $C_1(\Omega)$ designed the Poincar\'e Sobolev constant. Noticing that  
\begin{equation}
\label{it15}
e^{-\frac{m}{C_1(\Omega)}t}\frac{d}{dt}\bigg(e^{\frac{m}{C_1(\Omega)}t}|w|^2_2\bigg)=\frac{d}{dt}|w|^2_2+\frac{m}{C_1(\Omega)}|w|^2_2\leq|f|_q^{\frac{1}{\delta}}\bigg(\frac{2r}{m}\bigg)^{2r\sigma(r)}|w|^2_1\,m\,C^{\frac{2}{\delta}}.
\end{equation}
and integrating  (\ref{it15}) on  $[0,t)$ we get
  \begin{equation}
\label{it16}
|w(t)|^2_2\leq|w(0)|^2_2+\|f\|_{L^\infty(0,\infty,L^q(\Omega))}^{\frac{1}{\delta}}\bigg(\frac{2r}{m}\bigg)^{2r\sigma(r)}m\,C^{\frac{2}{\delta}}|w|^2_1.
\end{equation}
Since 
\begin{equation}
\label{it17}
|w(0)|^2_2=\int_\Omega w(0)^2dx=\int_\Omega u(0)^{2r}dx\leq|\Omega||u(0)|^{2r}_\infty\leq|\Omega|\tilde{U}_r^{2r},
\end{equation}
(\ref{it16}) and (\ref{it17}) gives us
\begin{equation}
\label{it18}
\tilde{U}_{2r}^{2r}\leq|\Omega|\tilde{U}_r^{2r}+\|f\|_{L^\infty(0,\infty,L^q(\Omega))}^{\frac{1}{\delta}}\bigg(\frac{2r}{m}\bigg)^{2r\sigma(r)}\,m\,C^{\frac{2}{\delta}}\tilde{U}_r^{2r}.
\end{equation}Whereas  $\frac{1}{\delta}>1$, $2r\sigma(r)>0$  and  $\sigma(r)=\frac{1}{2r\delta}$  it follows that 
\begin{equation}
\label{it20}
\tilde{U}_{2r}\leq C_2^{\sigma(r)}\|f\|_{L^\infty(0,\infty,L^q(\Omega))}^{\sigma(r)}r^{\sigma(r)}\tilde{U}_r,
\end{equation} with $C_2=C_2(\Omega,m)$. This completes the proof of Lemma.
\end{proof}

We have also
\begin{lemma}
\label{moser2}
Let $r>1$, $n\geq 3$, $p<\frac{n}{n-2}$ and $\sigma(r)=\frac{p(n+2)}{2[r(2p-pn+n)+np]}$  then we get
$$\sigma(2^kr)\leq\theta^k\sigma(r)\quad\forall k\in\mathbb{N},$$ with  $\theta\in(0,1).$
\end{lemma}
\begin{proof}By asking $c_1=\frac{p(n+2)}{2}$, $c_2=(2p-pn+n)$ and $c_3=np$ yields $\sigma(r)=\frac{c_1}{rc_2+c_3}$ with $c_1,c_2,c_3\in \mathbb{R}^{\star}_{+}$ .By taking  $\theta=1-\frac{c_2}{2c_2+c_3}$ the proof of this lemma is deduced by reasoning by induction.
\end{proof}

Returning now to the proof of the theorem, by lemma\,\ref{moser1} we get\begin{proof} 
$$\tilde U_{2r}\leq [C_2\,\|f\|_{L^\infty(0,\infty,L^q(\Omega))}]^{\sigma(r)}r^{\sigma(r)}\tilde U_{r}.$$By iterating this equation by taking $r=h,r=2h,r=2^2h,etc$, we obtain 
$$\tilde U_{2^{k+1}h}\leq [C_2\,\|f\|_{L^\infty(0,\infty,L^q(\Omega))}]^{\lambda1}\,2^{\lambda2}\,h^{\lambda1}\,\tilde U_{h},$$ with 
$$\lambda_1:=\sigma(h)+\sigma(2h)+\sigma(2^2h)+..+\sigma(2^{k-1}h)+\sigma(2^kr)$$et $$\lambda_2:=\sigma(2h)+2\sigma(2^2h)+3\sigma(2^3h)+...+(k-1)\sigma(2^{k-1}h)+k\sigma(2^kr).$$To complete the proof we just need to show that $\lambda_1,\lambda_2<+\infty.$ Indeed by lemma\,\ref{moser2}$$\lambda_1\leq\sum_{\mu=0}^k\alpha^{\mu}\sigma(h)\leq\sum_{\mu=0}^\infty\alpha^{\mu}\sigma(h)=\frac{\sigma(h)}{(1-\alpha)}<\infty.$$Noting also that $$\sigma(2^kh)\leq\theta^{k-1}\sigma(2h)\quad\forall k\in\mathbb{N}^\star,$$  it follows
$$\lambda_2\leq\sum_{\mu=1}^k\mu\alpha^{\mu-1}\sigma(2h)\leq\sum_{\mu=1}^\infty\mu\alpha^{\mu-1}\sigma(2h)=\frac{\sigma(2h)}{(1-\alpha)^2}<\infty.$$ This completes the proof of the theorem.
\end{proof}
\subsection{ Uniform estimate in time}
We prove in what follows an estimate of $u$ in $L^{\infty}(\mathbb{R}^+,H^1_0(\Omega)).$ We get
\begin{theorem}
\label{th1}
Assume that $f\in L^2(\Omega)$, $g\in H^1(\Omega)$, $u_0\in H^1_0(\Omega)$  and $a\in W^{1,\infty}(\mathbb{R})$ with $\displaystyle{\inf_{\mathbb{R}} a> 0}$. Then a solution $u$ of (\ref{eq0}) is such that  $u\in L^{\infty}(\mathbb{R}^+,H^1_0(\Omega)).$

\end{theorem}

 Taking a spectral basis related to the Laplace operator in the Galerkin approximation (see\,\cite{ref14}) we find that $-\Delta u$ can be regarded  as test function in  $L^2(0,T,L^2(\Omega))$ for all $T>0.$ By multiplying  (\ref{eq0}) by  $-\Delta u(t)$ and integrating over $\Omega$ 
\begin{equation}
\label{22oct1}
( u_t,-\Delta u)+( -div(a(l_r( u))\nabla u),-\Delta u)=(f,-\Delta u),
\end{equation} and also 
\begin{equation}
\label{23oct1}
\frac{1}{2}\frac{d}{dt}\|u\|^2_V+(-a(l_r(u))\Delta u ,-\Delta u)+(-a^{\prime}(l_r(u))\nabla l_r(u).\nabla u ,-\Delta u)=(f,-\Delta u).
\end{equation}Here $(.,.)$ is the usual scalar product on $L^2(\Omega)$. Taking to account 
\begin{equation}
|\nabla l_r(u)|_2\leq K\,\|g\|_{H^1(\Omega)}|\nabla u|_2,
\end{equation}where $K$ is a constant depending of $\Omega.$
It comes
\begin{equation}
\label{23oct3}
|(-a^{\prime}(l_r(u))\nabla l_r(u).\nabla u,-\Delta u)|\leq K\,\|g\|_{H^1(\Omega)}|a^{\prime}|_{\infty}\| u\|^2_V|\Delta u|_2
\end{equation} Now  from  (\ref{23oct3}) and (\ref{23oct1}) we have 
\begin{equation}
\label{23oct6}
\frac{1}{2}\frac{d}{dt}\|u\|^2_V+ m|\Delta u|_2^2-K\,\|g\|_{H^1(\Omega)}|a^{\prime}|_{\infty}\| u\|^2_V|\Delta u|_2\leq|f|_2|\Delta u|_2.
\end{equation}By using Young's inequality $ab\leq \frac{1}{2m}a^2+\frac{m}{2}b^2$, we get

\begin{equation}
\label{23oct8}
\frac{d}{dt}\|u\|^2_V\leq\frac{1}{m}|f|_2^2+\frac{1}{m} (K\,\|g\|_{H^1(\Omega)})^2\|a^{\prime}\|_{\infty}^2\| u\|^4.
\end{equation}In order to apply the uniform Gronwall lemma to (\ref{23oct8}) we start with a small estimate. Remember that 
\begin{equation}
\label{29oct3}
\frac{d}{dt}|u|^2_2+\,m\|u\|^2_V\leq\frac{1}{\lambda\,m}|f|_{2}^2,
\end{equation}where $\lambda$ is the principal eigenvalue of the Laplacian operator with Dirichlet boundary conditions.

 By integrating on $[t,t_0)$ we get 
\begin{equation}
|u(t+t_0)|^2_2+\,m\int_t^{t+t_0}\|u\|^2_V\,ds\leq\int_t^{t+t_0}\frac{1}{\lambda\,m}|f|_{2}^2\,ds+|u(t)|^2_2,
\end{equation}and also 
\begin{equation}
\label{bazo}
\int_t^{t+t_0}\|u\|^2_V\,ds\leq\frac{t_0}{\lambda\,m^2}|f|_{2}^2\,ds+\frac{1}{m}|u(t)|^2_2.
\end{equation}Let $\rho_0>0$ such that $|u(t)|^2_2\leq\rho_0^2 $. By setting 
\begin{equation*}
a_1=\frac{1}{m}c_1(\Omega)^2|a^{\prime}|_{\infty}^2a_3\quad a_2=\frac{t_0}{m}|f|_2^2\quad a_3=\frac{t_0\lambda}{m^2}|f|_{2}^2+\frac{1}{m}\rho_0^2,
\end{equation*}we obtain by using uniform Gronwall lemma to  (\ref{23oct8})
\begin{equation}
\label{29oct1}
\|u(t+t_0)\|_V\leq (\frac{a_3}{t_0}+a_2)exp(a_1)\quad\forall t\geq 0,\quad t_0>0.
\end{equation}Hence  $u\in L^{\infty}( t_0,+\infty,H^1_0(\Omega))$.
By using (\ref{23oct8}) and the classical Gronwall lemma it is easy to see that 
$u\in L^{\infty}(0,t_0,H^1_0(\Omega))$.  This completes the proof of the theorem.

\begin{remark}
\label{bab}
This theorem show us the existence of absorbing set in $H^1_0(\Omega)$.
 By considering $S(t)$ the semigroup associated to the equation\,(\ref{eq0})  defined by
\begin{equation*}
\begin{split}
S(t)\,:L^2(\Omega)&\rightarrow L^2(\Omega)\\
       u_0&\mapsto u(t), 
\end{split}
\end{equation*}with $u(t)$ a solution of (\ref{eq0}). As a result of the theorem\,\ref{th1} and the compact embedding of $H^1_0(\Omega)$ into $L^2(\Omega)$  we deduce that the semigroup $S(t)$  possesses a global attractor. Indeed it is easy to show the existence of absorbing set in $L^2(\Omega)$, the main difficulty here is to show that $S(t)$ is such that   
\begin{equation}
\begin{split}
&\forall B\subset L^2(\Omega)\quad\text{bounded},\quad\exists t_0=t_0(B)\\
&\text{such that}\quad\displaystyle{\bigcap_{t\geq t_0}\bigcup S(t)B}\quad\text{is relatively compact in}\quad L^2(\Omega).
 \end{split}
\end{equation}This property known as $S(t)$ is uniformly compact for $t$ large can be proved by using theorem\,\ref{th1} and the compact embedding of $H^1_0(\Omega)$ into $L^2(\Omega)$.
\end{remark}

\subsection{Asymptotic behaviour}
In this part we are interested in asymptotic behaviour of a weak solution of (\ref{eq0}). Our main interest here is the radial solutions. By radial solutions we means $\tilde{u}(t,|x|)=u(t,x)$. As in the stationary case $\Omega$ is a open ball of $\mathbb{R}^n$. Remember that \begin{equation*}
L^2_r(\Omega)=\{v\in L^2(\Omega)\quad\exists\tilde{v}\in L^2(]0,d/2[)\quad \text{such that}\quad v(x)=\tilde{v}(\|x\|) \}.
\end{equation*}In order to not make confusion between $u_0$ the solution to ($P_0$) and the initial value of (\ref{eq0}), we will take $u^0$  the initial value of (\ref{eq0}). 
\begin{theorem}
Assume that $f, g\in L^2_r(\Omega) $, $a$ is a continuous function and the assumption (\ref{p1}) checked then (\ref{eq0}) admits a radial solution.
\end{theorem}
\begin{proof}Let $w\in L^2(0,t, L^2_r(\Omega))$ it is clear that $l_r(w)$ is radial and also $a(l_r(w))$. Thus by (\ref{eq8}) $F_r$ maps $L^2(0,t, L^2_r(\Omega))$  into itself. The proof now follows by using arguments similar to those used  in theorem\,\ref{terib1}.
\end{proof}
Assume now 
\begin{equation}
\label{inf1}
f,g\geq 0\quad\text{in}\quad\Omega
\end{equation}and
 \begin{equation}
\label{inf2}
u_0\leq u^0\leq u_d,
\end{equation}with $u^0$ the initial value to (\ref{eq0}) and $u_0$ and $u_d$ respectively the solution of ($P_0$) and ($P_d$).

We can now give a stability result assuming that (\ref{eq0}) admits a unique solution.

\begin{theorem}Assume (\ref{inf1}) and $f, g\in L^2_r(\Omega)$.
Let $u$, $u_d$ and $u_0$ respectively the solution of (\ref{eq0}),  ($P_d$) and  ($P_d$). If 
\begin{equation*}
u_0\leq u^0\leq u_d,
\end{equation*} then 
\begin{equation*}
u_0\leq u\leq u_d\quad \forall t.
\end{equation*}
 \end{theorem}
\begin{proof}

Let 
\begin{equation}
\mathcal{S}=\{t\,|\quad l(u(s))\in [0,l_d(u_d)]\quad\forall s\leq t\}.
\end{equation}It is easy to prove that $\mathcal{S}$ contains 0 (see \ref{inf2}). By setting

\begin{equation}
t^{\star}=sup\{t\,|t\in\mathcal{S}\}.
\end{equation}
By continuity of the mapping  $t\,\rightarrow l_d(u(t))$, we get
\begin{equation}
 l_d(u(t^{\star}))\in [0,l_d(u_d)].
\end{equation}
By using (\ref{eq0}) and ($P_d$) we get in $\mathcal{D}(0,t^{\star})$
\begin{equation}
\label{inf4}
 \frac{d}{dt}(u_d-u,\phi)+\int_{\Omega}a(l_d(u))\nabla(u_d-u)\nabla\phi=-\int_{\Omega}(a(l_d(u_d))-a(l_d(u)))\nabla u_d\nabla\phi\quad\forall\phi\in H^1_0(\Omega).
\end{equation}
Choising $\phi=(u_d-u)^-$, (\ref{inf4}) becomes
\begin{equation}
\label{inf5}
 \frac{1}{2}\frac{d}{dt}|(u_d-u)^-|^2_2+\int_{\Omega}a(l_d(u_d))|\nabla(u_d-u)^-|^2=\int_{\Omega}(a(l_d(u))-a(l_d(u_d)))\nabla u_d\nabla(u_d-u)^-.
\end{equation}Since $a$ is nonincreasing ($a(l_d(u))-a(l_d(u_d))\geq 0\quad\forall t\leq t^{\star}$) hence proposition\,\ref{marche} yields 
\begin{equation}
\int_{\Omega}(a(l_d(u))-a(l_d(u_d)))\nabla u_d\nabla(u_d-u)^-\leq 0.
\end{equation}
Thus
\begin{equation}
\label{inf6}
\frac{1}{2}\frac{d}{dt}|(u_d-u)^-|^2_2+a(l_d(u_d))|\nabla(u_d-u)^-|^2_2\leq 0
\end{equation} Applying Poincarré Sobolev inequality we get
\begin{equation}
\frac{1}{2}\frac{d}{dt}|(u_d-u)^-|^2_2+C_2|(u_d-u)^-|^2_2\leq 0,  
\end{equation}this prove
\begin{equation}
\frac{d}{dt}\{e^{2t\,C_2}|(u_d-u)^-|^2_2\}\leq 0.  
\end{equation}
Moreover $(u_d-u)^-(0)=(u_d-u^0)^-=0$ it follows that $u_d\geq u\quad\forall t\in [0,t^{\star}]$.  In the same way we can also prove $u_0\leq u\quad\forall t\in [0,t^{\star}]$. It follows  
\begin{equation}
\label{inf9}
u_0\leq u\leq u_d\quad\forall t\in [0,t^{\star}]
\end{equation} To finish we just need to prove that $t^{\star}$ is very large,
this is typically the case. Indeed if $t^{\star}<\infty$ then 
\begin{equation}
l(u(t^{\star}))=0\quad\text{or}\quad l_d(u_d).
\end{equation}From (\ref{inf1}) and (\ref{inf9}) we deduce
\begin{equation}
\label{inf10}
u(t^{\star})= u_0\quad\text{or}\quad u(t^{\star})=u_d.
\end{equation}Due to the uniqueness of (\ref{eq0}), we deduce that $t=\infty$. This shows that 
\begin{equation*}
u_0\leq u\leq u_d\quad \forall t,
\end{equation*} and achieve the proof.
\end{proof}
\begin{remark}
The fact that $|u(t)|^2_2$ is not a Lyapunov function that is to say decreases in time, makes very complex the study of  certain asymptotic properties of our problem. Indeed under our study it is tempting to show that for $r$ fixed $r\in]0,d[$ \begin{equation*}
u(t)\rightarrow u^1_r\quad\text{in}\quad L^2(\Omega),
\end{equation*} where $u$ is the solution of (\ref{eq0}) and $u^1_r$ the solution belonging to the stable global branch described previously. A numerical study would be a great contribution to try to carry out some of our theoretical intuitions.
\end{remark}

\end{document}